\documentclass[12pt,a4paper]{article}

\usepackage[margin=1in]{geometry}  
\usepackage{graphicx}              
\usepackage{amsmath,euscript}               
\usepackage{amsfonts}              
\usepackage{amssymb,,amscd,epsf,amsthm, epsfig}                
\usepackage[toc,page]{appendix}

\newtheorem{thm}{Theorem}[section]

\newtheorem{prop}[thm]{Proposition}
\newtheorem{cor}[thm]{Corollary}
\newtheorem{lemma}[thm]{Lemma}
\newtheorem{claim}[thm]{Claim}
\theoremstyle{definition}
\newtheorem{rem}[thm]{Remark}
\newtheorem{defn}[thm]{Definition}

\input xy
\xyoption{all}

\newcommand{\Z}{\mathbb{Z}}
\newcommand{\LL}{\mathbb{L}}

\newcommand{\calF}{\mathcal{F}}
\newcommand{\calG}{\mathcal{G}}

\newcommand{\calK}{\mathcal{K}}

\newcommand{\calN}{\mathbb{N}}
\newcommand{\calP}{\mathcal{P}}

\newcommand{\DD}{\EuScript D}

\newcommand{\KK}{\EuScript K}

\newcommand{\pf}{{\it Proof.}\hspace{2ex}}

\newcommand{\epf}{\hspace*{\fill}\mbox{$\halmos$}}
\newcommand{\halmos}{\rule{1ex}{1.4ex}}

\newcommand{\proj}{\mathop{\mathrm{proj}}\nolimits}
\newcommand{\Hom}{\mathop{\mathrm{Hom}}\nolimits}

\newcommand{\Imm}{\mathop{\mathrm{Im}}\nolimits}

\newcommand{\modu}{\mathop{\mathrm{mod}}\nolimits}

\newcommand{\id}{\mathop{\mathrm{id}}\nolimits}
\newcommand{\inj}{\mathop{\mathrm{inj}}\nolimits}
\newcommand{\op}{\mathop{\mathrm{op}}\nolimits}
\newcommand\M          {{\mathcal{M}}}

\newcommand{\Tor}{\mathop{\mathrm{Tor}}\nolimits}

\newcommand{\Ext}{\mathop{\mathrm{Ext}}\nolimits}
\newcommand{\sg}{\mathop{\mathrm{sg}}\nolimits}

\newcommand{\HH}{\mathop{\mathrm{HH}}\nolimits}

\newcommand{\fGd}{\mathop{\mathrm{fGd}}\nolimits}

\begin{document}

\nocite{*}

\title{Singular Equivalence of Morita Type with Level}
\date{}

\author{Zhengfang WANG \thanks{
zhengfang.wang@imj-prg.fr, Universit\'e Paris Diderot-Paris 7, Institut
de Math\'ematiques de Jussieu-Paris Rive Gauche CNRS UMR 7586, B\^atiment Sophie Germain, Case 7012,
75205 Paris Cedex 13, France
}
}

\maketitle

\begin{abstract}
  We generalize the notion of stable equivalence of Morita type and define what is called ``singular equivalence of Morita type with level". Such an equivalence induces an equivalence between singular categories. We will also prove that a derived equivalence of standard type induces a
  singular equivalence of Morita type with level.
\end{abstract}

\section{Introduction}
Let $k$ be a commutative ring and let $A$ be a $k$-algebra. We
denote by $A$-$\modu$ the category of all finitely presented left $A$-modules,
and by $A$-$\underline{\modu}$ the stable module category of $A$-$\modu$ modulo projective
modules. The singular category $\DD_{\sg}(A)$ of $A$ is defined to be the Verdier quotient of the
bounded derived category $\DD^b(A)$ of finitely presented modules over $A$ by
the full subcategory $\KK^b(\mbox{$A$-$\proj$})$ consisting of bounded complexes of finitely presented projective
$A$-modules. Two $k$-algebras $A$ and $B$ are said to be stably equivalent if their stable categories
$A$-$\underline{\modu}$ and $B$-$\underline{\modu}$ are equivalent as $k$-categories,
and to be singularly equivalent if their singular categories $\DD_{\sg}(A)$ and $\DD_{\sg}(B)$
are equivalent as triangulated categories.

A stable equivalence of Morita type introduced by Brou\'e \cite{Brou} induces an equivalence of
stable categories in case $k$ is a field.
\begin{defn}[\cite{Brou}]
  Let $A$ and $B$ be two $k$-algebras. We say that $(_AM_B, _BN_A)$ defines a stable equivalence of Morita type if the following conditions are satisfied,
  \begin{enumerate}
    \item $M$ is finitely generated projective as a left $A$-module and as a right $B$-module.
    \item $N$ is finitely generated projective as a right $B$-module and as a left $A$-module.
    \item $M\otimes_BN\cong A\oplus P$ for some finitely generated projective $A$-$A$-bimodule $P$,
    and $N\otimes_AM\cong B\oplus Q$ for some finitely generated projective $B$-$B$-bimodule $Q$.
  \end{enumerate}
\end{defn}
Very recently analogous to the notion of stable equivalences
of Morita type, Xiao-Wu Chen and Long-Gang Sun defined in \cite{ChSu}
the concept of singular equivalences of Morita type.
\begin{defn}[\cite{ChSu}]\label{def-chen}
   Let $A$ and $B$ be two Noetherian $k$-algebras. We say that $(_AM_B, _BN_A)$ defines a singular equivalence of Morita type if the following conditions are satisfied,
  \begin{enumerate}
    \item $M$ is finitely generated projective as a left $A$-module and as a right $B$-module.
    \item $N$ is finitely generated projective as a left $B$-module and as a right $A$-module.
    \item There are bimodule isomorphisms $M\otimes_BN\cong A\oplus P, $ where $P$ is finitely presented and of finite projective dimension as an $A$-$A$-bimodule
    and $N\otimes_AM\cong B\oplus Q,$ where $Q$ is finitely presented and of finite projective dimension as a $B$-$B$-bimodule.
  \end{enumerate}
\end{defn}
If $k$ is a field and $(M, N)$ defines a singular equivalence of Morita type then $$M\otimes_B-: \DD_{\sg}(B)\rightarrow \DD_{\sg}(A)$$
is an equivalence of triangulated categories with quasi-inverse
$$N\otimes_A-:\DD_{\sg}(A)\rightarrow \DD_{\sg}(A).$$

In this paper, we generalize Chen and Sun's notion of singular equivalence of Morita type and define a singular equivalence of Morita type with level (cf. Definition \ref{defn-Main-def} below).
This new concept is very related to
derived equivalences, that is, a derived equivalence of standard type induces a singular
equivalence of Morita type with level. This generalizes the fact in \cite{Ric} and \cite{KeVo} that derived equivalences between two self-injective $k$-algebras induce stable
equivalences of Morita type.

Much work has be done to study the invariants under stable equivalence of Morita type
(eg. \cite{Liu}, \cite{LiXi}, \cite{PaZh}, \cite{Pogo})
and singular equivalence of Morita type (eg. \cite{ChSu}, \cite{ZhZi}). Singular equivalences of Morita type with level also preserve some invariants, for example, the Hochschild homology group
in positive degree (cf. Proposition \ref{prop-Hoschild-Homology}) and the Gorensteinness of algebras. In general, self-injectivity of algebras is
not preserved under singular equivalences of Morita type with level (cf. Remark \ref{rem-expl-self}), but it is preserved under singular equivalences of
Morita type (in the sense of Chen and Sun) (cf. Proposition \ref{prop-singular-Sun}).

The paper is organised as follows. In Section 2, we define the concept of a singular equivalence of Morita type
with level and show that a derived equivalence of standard type induces a singular equivalence of Morita type with level.
An example of a singular equivalence of Morita type with level is given in Section 3.
This example illustrates that Hochschild cohomology in lower degree
may not be invariant under singular equivalences of Morita type with level. Section 4 is devoted to study the category $\underline{\calG_A}$ which is the image in $A$-$\underline{\modu}$ of the category $\calG_A$ consisting of modules with Gorenstein dimension zero (Definition \ref{defn-Goren-zero}). We show that a singular equivalence of Morita type with level between two $k$-algebras $A$ and $B$ induces an equivalence between $\underline{\calG_A}$ and $\underline{\calG_B}$ under some conditions. Moreover,  we will give a sufficient condition about when a singular
equivalence of Morita type with level induces a stable equivalence of Morita type.
In Appendix A, we show that a singular equivalence of Morita type (defined by Chen and Sun) preserves self-injectivity of algebras, but in general
a singular equivalence of Morita type with level does not.

In this paper, we need to use theories about derived categories and
two-sided tilting complexes, for details, we refer to \cite{Kel}, \cite{Ric} and \cite{Zim}.
For general homological theory, we refer to \cite{Wei} and \cite{Zim}.
Throughout this paper, $k$ denotes a commutative ring (or sometimes a field when explicitly specified).
$A$ and $B$ denote $k$-algebras which are finite generated and projective as $k$-modules.

\bigskip
\noindent
{\bf Acknowledgements:} This work is a part of author's PhD thesis, I would like to thank my PhD supervisor Alexander Zimmermann for introducing this interesting topic and for his many valuable suggestions for improvement. I want to thank Guodong Zhou for many inspiring conversations and
constant supports. I also want to thank Yuming Liu for many discussions and Hirotaka Koga for
an interesting discussion in the 16th International Conference of Representation of Algebras and Workshop (ICRA 2014)
in Sanya.

\section{Singular equivalence of Morita type with level}
It is well known that an algebra $A$ is of finite global dimension if and only if $\DD_{\sg}(A)=0$ (cf. \cite{Hap},
\cite{ZhZi}, \cite{Zim},). In this paper,  all the algebras we consider are assumed to be of infinite global dimension except where noted otherwise. For simplicity, we denote $\otimes_k$ by $\otimes$, and denote the enveloping algebra $A\otimes A^{\op}$ of $A$ by $A^e$.
\begin{defn}\label{defn-Main-def}
  Let $k$ be a commutative algebra and let $A$ and $B$ be two $k$-algebras which are projective as $k$-modules. Let $M$ and $N$ be, respectively, an $A$-$B$-bimodule and a $B$-$A$-bimodule.
  We say that $(M, N)$ defines a singular equivalence of Morita type with level $n$ for some $n\in \calN$ if the following conditions are satisfied:
  \begin{enumerate}
    \item $M$ is projective as a left $A$-module and as a right $B$-module,
    \item $N$ is projective as a left $B$-module and as a right $A$-module,
    \item $M\otimes_B N\cong \Omega^n_{A^e}(A)$ in $A^e$-$\underline{\modu}$ and $N\otimes_AM\cong \Omega^n_{B^e}(B)$ in $B^e$-$\underline{\modu}$.
  \end{enumerate}
\end{defn}
\begin{rem}
  We remark that if $k$ is a field, then the tensor functor $$M\otimes_B-: \DD_{\sg}(B)\rightarrow \DD_{\sg}(A)$$
  is an equivalence with quasi-inverse $$N\otimes_A-:\DD_{\sg}(A)\rightarrow \DD_{\sg}(B).$$
  Moreover, $n=0$ gives the case of a stable equivalence of Morita type.
\end{rem}
Next we prove that a derived equivalence of standard type induces a singular equivalence of Morita type with level.
\begin{thm}\label{thm1}
  Let $k$ be a commutative ring. Let $A$ and $B$ be two Noetherian $k$-algebras of infinite global dimension such that $A$ and $B$ are derived equivalent. Suppose that $A$ and $B$ are projective as $k$-modules.  Then there exist an $A$-$B$-bimodule $M$ and a $B$-$A$-bimodule $N$ such that $(M, N)$ defines a singular equivalence of Morita type with level. More precisely, an equivalence of standard type induces a singular equivalence of Morita type with level $n$  for some $n$ which is larger than the length of the tilting complex inducing the derived equivalence.
\end{thm}
\pf Since $\DD^b(A)\cong\DD^b(B)$, there exist two-sided tilting complexes $X\in \DD^b(A\otimes B^{\op})$ and $Y\in \DD^b(B\otimes A^{\op})$ such that (cf. \cite{Kel1}, \cite{Ric})
\begin{eqnarray*}
X\otimes_B^{\mathbb{L}}Y&\cong&A
\end{eqnarray*}
in $\DD^b(A\otimes A^{\op})$ and
\begin{eqnarray*}
Y\otimes_A^{\mathbb{L}}X&\cong &B
\end{eqnarray*}
in $\DD^b(B\otimes B^{\op})$.
Moreover, from [Proposition 6.2.27, \cite{Zim}] we can choose $X$ of the following form:
$$0\leftarrow X_k\leftarrow X_{k+1}\leftarrow \cdots \leftarrow X_{l} \leftarrow0$$
where
$X_k, \cdots, X_{l-1} $ are projective as $A$-$B$-bimodules and the most right term $X_l$ is projective as a left $A$-module and as a right $B$-module.
Similarly, choose $Y$ of the form:
$$ 0\leftarrow Y_p\leftarrow Y_{p+1}\leftarrow \cdots \leftarrow Y_{q}\leftarrow 0,$$
where
$Y_p, \cdots, Y_{q-1} $ are projective as $B$-$A$-bimodules and the most right term $Y_q$ is projective as a right $A$-module and as a left $B$-module.

Then $$X\otimes^{\mathbb{L}}_BY\cong  (0\leftarrow Z_{k+p}\leftarrow \cdots \leftarrow Z_{l+q}\leftarrow 0)$$
where $Z_{m}:=\bigoplus_{i+j=m} X_{i}\otimes_B Y_{j}$.
\begin{claim}
We claim that $Z_{m}$ is projective as $A$-$A$-bimodules for all $m<q+l$.
\end{claim}
It is sufficient to prove that $X_{i}\otimes_B Y_{j}$ is projective for $i+j<q+l$.
Since $i+j<q+l$, we know either $X_{i}$ or $Y_{j}$ is projective as bimodules. Without loss of generality, we assume that $X_{i}$ is projective as $A$-$B$-bimodule, thus $X_{i}\otimes_BY_{j}$ is projective as $A$-$A$-bimodule since $Y_j$ is projective as  a right $A$-module. Therefore $Z_m$ is
projective as an $A$-$A$-bimodule for $m<l+q$.

Since $X\otimes_B^{\mathbb{L}}Y\cong A$ in $\DD^b(A\otimes A^{\op})$, the homology group
 $$H^i(X\otimes^{\LL}_BY)=
 \begin{cases}
0 & i\neq 0\\
A & i=0,
 \end{cases}
 $$
 hence we have the following exact sequences:
\begin{eqnarray} \label{exact-seq1}
0\leftarrow Z_{k+p}\leftarrow \cdots\leftarrow Z_{-1}\leftarrow Z_0\leftarrow \ker d_0\leftarrow 0,
\end{eqnarray}
\begin{eqnarray}\label{exact-seq2}
0\leftarrow A\cong \ker d_0/\Imm d_1\leftarrow \ker d_0\leftarrow Z_1\leftarrow \cdots\leftarrow Z_{l+q}\leftarrow 0.
\end{eqnarray}

From the exact sequence (\ref{exact-seq1}) and the fact that the modules $Z_{m}$ are projective as bimodules for $k+p \leqslant m<l+q$, it follows that the exact sequence (\ref{exact-seq1}) splits, thus $\ker d_0$ is projective as bimodule. Therefore from the exact sequence (\ref{exact-seq2}) and the assumption that $A$ is not of finite global dimension we obtain that in $A^e$-$\underline{\modu}$ $$\Omega_{A^e}^{l+q+1}(A)\cong Z_{l+q}\cong X_{l}\otimes_B Y_{q}.$$
Similarly, we can get that in $B^e$-$\underline{\modu}$ $$\Omega_{B^e}^{l+q+1}(B)\cong Y_{q}\otimes_A X_{l}.$$
Hence, $(X_l, Y_q)$ defines a singular equivalence of Morita type with level $l+q+1$.
\epf

\section{Examples}
In this section, we will give an example of a singular equivalence of Morita type with level.

Let $k$ be a field. Let $A$ be a finite-dimensional $k$-algebra of infinite global dimension and $S$ be a finite-dimensional $k$-algebra which has finite projective dimension $m_1$ as $S$-$S$-bimodule, that is, $\Omega_{S\otimes S^{\op}}^{m_1-1}(S)$ is not projective  and $\Omega_{S\otimes S^{\op}}^{m_1}(S)$ is projective as bimodules. Let $M$
be an $A$-$S$-bimodule with finite projective dimension $m_2$ as a bimodule. Set $n:=\max\{m_1, m_2\}$,  $$B:=\left(\begin{smallmatrix}
  S &0 \\ M &A
\end{smallmatrix}\right)
$$
and $$\epsilon:=\left(\begin{smallmatrix}
  0 & 0\\0&1
\end{smallmatrix}\right).$$
We claim that $$(\Omega_{A\otimes B^{\op}}^{n+1}(\epsilon B), B\epsilon)$$ defines a
singular equivalence of Morita type with level $n+1$ between $A$ and $B$. The following lemmas are devoted to prove this claim.

\begin{lemma}\label{lemma-fully}
  $\Omega_{A\otimes B^{\op}}^{n+1}(\epsilon B)$ is projective as a left $A$-module and
  as a right $B$-module. $B\epsilon$ is projective as a left $B$-module and as a right $A$-module.
\end{lemma}
\pf Observe that $\epsilon B \epsilon \cong A$ as algebras, hence $\epsilon B$ is an $A$-$B$-bimodule. Since $\epsilon B=\left( \begin{smallmatrix} 0 & 0\\ M &A \end{smallmatrix}\right)$, $\epsilon B\cong M \oplus A$ as $A$-modules, hence we have the following isomorphism in
$A$-$\underline{\modu}$,
\begin{eqnarray*}
  \Omega_{A\otimes B^{\op}}^{n+1}(\epsilon B)&\cong& \Omega_A^{n+1}(\epsilon B)\\
  &\cong& \Omega_A^{n+1}(M)\\
  &\cong & \Omega_{A\otimes B^{\op}}^{n+1}(M)\\
  &\cong & 0
\end{eqnarray*}
where the first and the third isomorphisms hold since projective $A\otimes B^{\op}$-modules are projective $A$-modules and syzygies  are independent of the projective resolutions.
Hence $\Omega_{A\otimes B^{\op}}^{n+1}(\epsilon B)$ is projective as a left $A$-module. As right $B$-modules we have the following isomorphism
$$\epsilon B\oplus (1-\epsilon)B \cong B,$$ hence $\epsilon B$ is a projective right $B$-module. So $\Omega_{A\otimes B^{\op}}^{n+1}(\epsilon B)$ is projective as a left $A$-module and as a right $B$-module. Since $$B\epsilon \oplus B(1-\epsilon)\cong B$$ as left $B$-modules, $B\epsilon $ is projective as a left $B$-module. As right $A$-modules we have $$B\epsilon \cong \epsilon B\epsilon \oplus (1-\epsilon) B\epsilon \cong \epsilon B\epsilon\cong A$$ since $(1-\epsilon) B\epsilon=0$.
Hence $B\epsilon $ is projective as a left $B$-module and as a right $A$-module.
\epf

\begin{rem}
  In fact, we have $$\epsilon B\otimes_B B\epsilon \cong A$$ as $A$-$A$-bimodules and
  $$B\epsilon\otimes_A\epsilon B \cong B\epsilon B$$ as $B$-$B$-bimodules.
\end{rem}

We define two functors:
\begin{eqnarray*}
\begin{tabular}{ccccc}
$\mathcal{F}$ : &\mbox{$A\otimes S^{\op}$-$\modu$}&$\rightarrow$ & \mbox{$B^e$-$\modu$}\\
 &$X$&$\mapsto$&$\mbox{$\left( \begin{smallmatrix}  0 & 0 \\ X & 0        \end{smallmatrix}   \right) $}$\\
\end{tabular}
\end{eqnarray*}
and
\begin{eqnarray*}
  \begin{tabular}{ccccc}
$\mathcal{G}$ : &\mbox{$S\otimes S^{\op}$-$\modu$}&$\rightarrow$ & \mbox{$B^e$-$\modu$}\\
 &$Y$&$\mapsto$&$\mbox{$\left( \begin{smallmatrix}  Y & 0 \\M\otimes_S Y & 0        \end{smallmatrix}   \right) $}$\\
\end{tabular}
\end{eqnarray*}

\begin{lemma}\label{lemma-full}
   Both $\mathcal{F}$ and $\mathcal{G}$ are fully-faithful and send projective objects to projective objects. Moreover $\mathcal{F}$ is exact.
\end{lemma}

\pf Observe that $$\calF\cong B\epsilon\otimes_A-\otimes_S(1-\epsilon)B.$$ It is clear that $B\epsilon\otimes_A-\otimes_S(1-\epsilon)B$ is exact and sends projective objects to
projective objects since $B\epsilon$ is projective as a left $B$-module and as a right $A$-module (cf. Lemma \ref{lemma-fully}) and $(1-\epsilon)B$ is projective as a left $S$-module and
as a right $B$-module.
It is clear that we have the following isomorphism for any $X_1, X_2\in A\otimes S^{\op}$,
\begin{eqnarray*}
\Hom_{A\otimes S^{\op}}(X_1, X_2)\cong \Hom_{B^e}(\left( \begin{smallmatrix}  0 & 0 \\ X _1& 0        \end{smallmatrix}   \right),  \left( \begin{smallmatrix}  0 & 0 \\ X_2 & 0        \end{smallmatrix}   \right))
\end{eqnarray*}
Therefore $\calF$ is fully-faithful, exact and sends projective objects to projective objects.
Similarly, we also observe that $$\calG\cong B(1-\epsilon)\otimes_S -\otimes_S (1-\epsilon)B,$$ hence $\calG$
preserves projective objects.
We define the functor
\begin{eqnarray*}
  \begin{tabular}{ccccc}
$\mathcal{K}$ : &\mbox{$B^e$-$\modu$}&$\rightarrow$ & \mbox{$S^e$-$\modu$}\\
 &$Y$&$\mapsto$& $S\otimes_BY\otimes_BS$,
\end{tabular}
\end{eqnarray*}
and we have $$\calK\circ \calG\cong \id.$$
Therefore $\calG$ is a faithful functor. On the other hand,
by computation, it is clear that $$
\Hom_{S^e}(X, Y)\rightarrow \Hom_{B^e}(\calG(X), \calG(Y))$$ is an epimorphism.
Therefore, $\calG$ is fully-faithful and sends projective objects to projective objects.
\epf

\begin{lemma}
  $\Omega_{B^{e}}^{n+1}(B)\cong \Omega_{B^e}^{n+1}(B\epsilon B)$ in $B^e$-$\underline{\modu}$.
\end{lemma}
\pf We have the following commutative diagram in $B^e$-$\modu$ (cf. Section 5 in \cite{GMS})
\begin{eqnarray*}
\xymatrix@C=1pc{
&&0& 0\\
0\ar[r]  &  \mbox{$\left( \begin{smallmatrix}  0 & 0 \\ M & 0        \end{smallmatrix}   \right) $} \ar[r] &\mbox{$\left( \begin{smallmatrix}  0 & 0 \\ M & A      \end{smallmatrix}   \right) $}\oplus
\mbox{$\left( \begin{smallmatrix}  S &0 \\ M & 0        \end{smallmatrix}   \right) $}  \ar[r] \ar[u] & \mbox{$\left( \begin{smallmatrix}  S& 0 \\ M & A       \end{smallmatrix}   \right) $} =:B\ar[r] \ar[u]& 0\\
&0\ar[r]& \mbox{$\left( \begin{smallmatrix}  0 & 0\\ 0 & A        \end{smallmatrix}   \right) $}\otimes  \mbox{$\left( \begin{smallmatrix}  0 & 0 \\ M & A       \end{smallmatrix}   \right) $} \oplus \mbox{$\left( \begin{smallmatrix}  S & 0\\ M & 0       \end{smallmatrix}   \right) $}\otimes \mbox{$\left( \begin{smallmatrix}  S & 0\\ 0 & 0       \end{smallmatrix}   \right) $}
\ar[u]\oplus K\ar[r]^-{\cong}&    B\otimes B\ar[u]^-{\mu_B}  \\
&0\ar[r]& \Omega^1_{B^e}(\mbox{$\left( \begin{smallmatrix}  0 & 0 \\ M & A      \end{smallmatrix}   \right) $})\oplus \Omega^1_{B^e}(\mbox{$\left( \begin{smallmatrix}  S & 0 \\ M & 0      \end{smallmatrix}   \right) $})\oplus K\ar[u] \ar[r] & \Omega_{B^e}^1(B)\ar[u]\\
&&0\ar[u]&0\ar[u]
}
\end{eqnarray*}
where $$K=\mbox{$\left( \begin{smallmatrix}  0 & 0\\ 0 & A        \end{smallmatrix}   \right) $}\otimes \mbox{$\left( \begin{smallmatrix}  S & 0\\ 0 & 0       \end{smallmatrix}    \right)$}\oplus \mbox{$\left( \begin{smallmatrix}  S & 0\\ M& 0       \end{smallmatrix}   \right) $}\otimes  \mbox{$\left( \begin{smallmatrix}  0 & 0 \\ M& A    \end{smallmatrix}   \right) $},$$ the top row is exact, the middle row is a decomposition of $B\otimes B^{\op}$ as a $B$-$B$-bimodule and  we denote, respectively,
the kernels of the natural morphisms
$$\mbox{$\left( \begin{smallmatrix}  0 & 0\\ 0 & A        \end{smallmatrix}   \right) $}\otimes  \mbox{$\left( \begin{smallmatrix}  0 & 0 \\ M & A       \end{smallmatrix}   \right) $}\rightarrow \mbox{$\left( \begin{smallmatrix}  0 & 0\\ M & A        \end{smallmatrix}   \right) $}$$ and $$ \mbox{$\left( \begin{smallmatrix}  S & 0\\ M & 0       \end{smallmatrix}   \right) $}\otimes \mbox{$\left( \begin{smallmatrix}  S & 0\\ 0 & 0       \end{smallmatrix}   \right) $}
\rightarrow \mbox{$\left( \begin{smallmatrix}  S & 0\\ M & 0        \end{smallmatrix}   \right) $}$$ induced by the multiplication $\mu_B$ by $\Omega^1_{B^e}(\mbox{$\left( \begin{smallmatrix}  S & 0 \\ M & 0     \end{smallmatrix}   \right) $})$ and $\Omega^1_{B^e}(\mbox{$\left( \begin{smallmatrix}  S & 0 \\ M & 0      \end{smallmatrix}   \right) $}).$

Therefore from the Snake Lemma, we have the following exact sequence in $B^e$-$\modu$
\begin{eqnarray*}
  0\rightarrow \Omega^1_{B^e}(\left( \begin{smallmatrix} 0 & 0\\ M &A \end{smallmatrix}\right)) \oplus \Omega^1_{B^e}(\mbox{$\left( \begin{smallmatrix}  S & 0 \\ M & 0      \end{smallmatrix}   \right) $})  \oplus K\rightarrow \Omega^1_{B^e}(B) \rightarrow \left( \begin{smallmatrix} 0 & 0\\ M &0 \end{smallmatrix}\right)\rightarrow 0.
\end{eqnarray*}
Hence we have the following distinguished triangle in the left triangulated category $B^e$-$\underline{\modu}$ since $K$ is projective as a $B$-$B$-bimodule
\begin{eqnarray*}
\Omega^1_{B^e}(\left( \begin{smallmatrix} 0 & 0\\ M &0 \end{smallmatrix}\right))\rightarrow \Omega^1_{B^e}(\left( \begin{smallmatrix} 0 & 0\\ M &A \end{smallmatrix}\right)) \oplus \Omega^1_{B^e}(\mbox{$\left( \begin{smallmatrix}  S & 0 \\ M & 0      \end{smallmatrix}   \right) $}) \rightarrow \Omega_{B^e}^1(B)\rightarrow \left( \begin{smallmatrix} 0 & 0\\ M &0 \end{smallmatrix}\right).
\end{eqnarray*}
 In left triangulated categories one may apply shift functors to the left on distinguished triangles, and hence we
 have the following distinguished triangle,
 \begin{eqnarray}\label{eqnarray-left}
\Omega^{n+1}_{B^e}(\left( \begin{smallmatrix} 0 & 0\\ M &0 \end{smallmatrix}\right))\rightarrow \Omega^{n+1}_{B^e}(\left( \begin{smallmatrix} 0 & 0\\ M &A \end{smallmatrix}\right)) \oplus \Omega^{n+1}_{B^e}(\mbox{$\left( \begin{smallmatrix}  S & 0 \\ M & 0      \end{smallmatrix}   \right) $}) \rightarrow \Omega_{B^e}^{n+1}(B)\rightarrow \Omega_{B^e}^n( \left( \begin{smallmatrix} 0 & 0\\ M &0 \end{smallmatrix}\right))
\end{eqnarray}

Now we use the functors $\mathcal{F}$ and $\mathcal{G}$ and we know that $$
\mathcal{F}(M)\cong \mbox{$\left( \begin{smallmatrix}  0 & 0 \\ M & 0      \end{smallmatrix}   \right) $}$$
$$\mathcal{G}(S)\cong \mbox{$\left( \begin{smallmatrix}  S & 0 \\ M & 0      \end{smallmatrix}   \right) $}.$$
Since $M$ and $S$ have finite projective dimension respectively, in $A\otimes S^{\op}$-$\modu$
and $S^e$-$\modu$, it follows from Lemma \ref{lemma-full}  that $\mathcal{F}(M)$ and $\mathcal{G}(S)$ both have finite
projective dimension smaller than $n+1$  in $B^e$-$\modu$, hence we get in $B^e$-$\underline{\modu}$,  $$\Omega^{n+1}_{B^e}(\mbox{$\left( \begin{smallmatrix}  S & 0 \\ M & 0      \end{smallmatrix}   \right) $})\cong 0$$
and $$\Omega_{B^e}^n( \left( \begin{smallmatrix} 0 & 0\\ M &0 \end{smallmatrix}\right))\cong 0.$$
So from the distinguished triangle (\ref{eqnarray-left}) above, we have the following isomorphism in $B^e$-$\underline{\modu}$
$$\Omega_{B^{e}}^{n+1}(B)\cong \Omega_{B^e}^{n+1}(B\epsilon B).$$
Therefore we have shown that  $$(\Omega_{A\otimes B^{\op}}^{n+1}(\epsilon B), B\epsilon)$$ defines a
singular equivalence of Morita type with level $n+1$ between $A$ and $B$.
\epf

\begin{rem}\label{rem-expl-self}
  This example above illustrates that the self-injectivity of algebras may not be preserved under singular equivalence of Morita type with level. Assume that $A$ is a self-injective $k$-algebra, in general, $B=\left(\begin{smallmatrix}
  S &0 \\ M &A
\end{smallmatrix}\right)
$ is not a self-injective algebra.
\end{rem}

\begin{rem}
We also remark that the example above illustrates that Hochschild cohomology of
lower degree may not be invariant under singular equivalence of Morita type
with level.

Let $k$ be a field. Let $A$ be a finite dimensional $k$-algebra and $S=k$, then $$B=\left(\begin{smallmatrix}
  k &0 \\ M &A
\end{smallmatrix}\right)$$
We know that there is a singular equivalence of Morita type with level (n+1) between $A$
and $B$.
From \cite{GMS}, there is an exact sequence:
$$\cdots\rightarrow \Ext^0(M, M)\rightarrow \HH^1(B)\rightarrow \HH^1(A)\rightarrow
\Ext^1(M, M)\rightarrow\HH^2(B)\rightarrow  \cdots$$
In general, $\Ext^i(M, M)\neq0$ for $i<n$, hence $\HH^i(A)$ is in general not isomorphic to $\HH^i(B)$ for $i<n$.

\end{rem}

Next we show that the Hochschild homology of positive degree is invariant under singular equivalence of Morita type with level.
\begin{prop}\label{prop-Hoschild-Homology}
Let $k$ be a field. Let $A$ and $B$ be two finite dimensional $k$-algebras which are singularly equivalent of Morita type with level.
Then $\HH_i(A)\cong \HH_i(B)$ for $i>0$.
\end{prop}
\pf The proof is similar to the proof of Theorem 4.1 in \cite{ZhZi}.
Recall that for an $A$-$B$-bimodule $M$, we can define a transfer map $t_M: \HH_i(A)\rightarrow \HH_i(B)$ for
each $i>0$. Note that $$(-1)^it_{\Omega_{A^e}^i(X)}=t_{X}$$ for any $A$-$A$-bimodule $X$ and $i\in\calN$.
Hence $(-1)^it_{\Omega_{A^e}^i(A)}=t_{A}=\id$. Then we have
$$t_M\circ t_N\cong (-1)^n\id$$ and $$t_N\circ t_M\cong (-1)^n \id.$$
So $t_M$ induce an isomorphism between $\HH_i(A)$ and $\HH_i(B)$ for $i>0$.
\epf

\section{Modules of Gorenstein dimension zero}
\subsection{Gorenstein dimension zero}
We assume that $k$ is a field throughout this subsection.
\begin{defn}\label{defn-Goren-zero}
Let $A$ be a finite-dimensional $k$-algebra and $M$ be a finitely generated $A$-module.  We say that $M$ has Gorenstein dimension zero if
it satisfies the following conditions,
\begin{enumerate}
\item The natural homomorphism $M\rightarrow \Hom_A(\Hom_A(M, A), A)$ is an isomorphism,
\item $\Ext_A^i(M, A)=\Ext_A^i(\Hom_A(M, A), A)=0$ for $i>0$.
\end{enumerate}
We denote by $\calG_A$ the full subcategory of $A$-$\modu$ consisting of modules of Gorenstein dimension zero.
Observe that $A$-$\proj\subset \calG_A$,  we denote the full subcategory consisting of the image of $\calG_A$ in $A$-$\underline{\modu}$ by $\underline{\calG_A}$ and we denote by $\DD^b(A)_{\fGd}$ the full subcategory consisting of  the image of $\KK^b(\calG_A)$ in $\DD^b(A)$. Clearly, $\KK^b(\mbox{$A$-$\proj$})\subset \DD^b(A)_{\fGd}.$
\end{defn}
\begin{rem}\label{rem-zero1}
We remark the fact that $\underline{\calG_A}$ is a triangulated category with the suspension functor $\Phi_A$, which is a quasi-inverse of $\Omega_A$ (cf. \cite{HoKo}, \cite{Avr}).
An $A$-module $M$ is of Gorenstein dimension zero if and only if there exists an acyclic complex of projective $A$-modules
\begin{eqnarray*}
  \xymatrix{P_{\bullet}: \cdots& P_{-1} \ar[l]_-{d_{-1}}& P_0\ar[l]_-{d_{0}} &P_1\ar[l]_-{d_1} &  P_2\ar[l]_-{d_2} & \cdots\ar[l]}
\end{eqnarray*}
such that $\ker(d_0)\cong M$ and $\Hom_A(P_{\bullet}, A)$ is acyclic.
\end{rem}
\begin{thm}[Theorem 4.2. \cite{Ka}]\label{thm-Ka1}
  Let $k$ be a field. Let $A$ and $B$ be two finite dimensional $k$-algebras. Then the equivalence
  $X\otimes^{\LL}_B-: \DD^b(B)\rightarrow \DD^b(A)$ induces an equivalence
  $$X\otimes^{\LL}_B-: \DD^b(A)_{\fGd}\rightarrow \DD^b(B)_{\fGd}.$$
\end{thm}

Next we will prove that under some condition,  a singular equivalence of Morita type with level between two $k$-algebras $A$ and $B$ induces an equivalence between $\underline{\calG_A}$ and
$\underline{\calG_B}$ (cf. Proposition \ref{prop-zero1}).

\begin{lemma}\label{lemma-defect1}
Let $A$ be a finite-dimensional $k$-algebra, $X\in \mathcal{G}_A$
and $M$ be any $A$-module such that $\Ext_A^i(X, M)=0$ for $i>0$. Then we have the following natural isomorphism $$\Hom_A(X, A)\otimes_AM\cong \Hom_A(X, M).$$
\end{lemma}
\pf Since $X\in \mathcal{G}$, we can take a right resolution of $X$, $$0\rightarrow X\rightarrow P_0\rightarrow P_1\rightarrow \cdots$$
Since $\Ext_A^i(X, M)=0$ for $i>0$, we have the following exact sequence
$$0\leftarrow \Hom_A(X, M)\leftarrow\Hom_A(P_0, M)\leftarrow \Hom_A(P_1, M)\leftarrow\cdots$$
Since $P_i$ are projective $A$-modules, $$\Hom_A(P_i, M)\cong \Hom_A(P_i, A)\otimes_AM.$$
Since $\Ext^i_A(X, A)=0$ for $i>0$, we have the following exact sequence
$$0\leftarrow \Hom_A(X, A)\leftarrow\Hom_A(P_0, A)\leftarrow \Hom_A(P_1, A)\leftarrow\cdots$$
Then we have the following commutative diagram
\begin{eqnarray*}
\xymatrix{
0 & \Hom_A(X, M) \ar@{->>}[l] & \Hom_A(P_0, M)  \ar[l] & \Hom_A(P_1, M)\ar[l] & \cdots\ar[l]\\
0&\Hom_A(X, A)\otimes_AM \ar@{->>}[l] \ar[u] &  \Hom_A(P_0, A)\otimes_A M \ar[l]\ar[u]^{\cong}  & \Hom_A(P_1, A)\otimes_A M  \ar[l]\ar[u]^{\cong} & \cdots\ar[l]
}
\end{eqnarray*}
where we have seen that the first row sequence is exact.
Hence we have $$\Hom_A(X, M)\cong\Hom_A(X, A)\otimes_AM.$$
\epf

\begin{prop}\label{prop-zero1}
Let $(_AM_B, _BN_A)$ define a singular equivalence of Morita type of level $n$ between $A$ and $B$ such that $\Hom_A(M, A)$ and $\Hom_B(N, B)$ are of finite projective dimension as a $B$-module and as an $A$-module respectively. Then we have the following two equivalences, which are quasi-inverse to each other,
$$M\otimes_B-: \underline{\mathcal{G}_B}\rightarrow \underline{\mathcal{G}_A}$$
and $$N\otimes_A-:\underline{\mathcal{G}_A}\rightarrow \underline{\mathcal{G}_B}$$
\end{prop}
\pf Let $X\in\mathcal{G}_B$, we need to show that $M\otimes_BX\in \mathcal{G}_A$, that is to say,
we need to prove the following:
\begin{enumerate}
\item $\Ext_A^i(M\otimes_BX, A)=0$ for $i>0$. Since $M$ is projective as a left $A$-module and as a right $B$-module,
 $$\Ext_A^i(M\otimes_B X, A)\cong \Ext_B^i(X, \Hom_A(M, A))$$
Since $\Hom_A(M, A)$ is of finite projective dimension as $B$-module and $X\in\mathcal{G}_B$, $$ \Ext_B^i(X, \Hom_A(M, A))=0,$$ where we need to use the induction on the projective dimension of $Y$ to argue that $\Ext_B^i(X, Y)=0$ for any $i>0$ and any $B$-module $Y$ of finite projective dimension since $X$ is of Gorenstein dimension zero.
Hence $$\Ext_A^i(M\otimes_B X, A)=0$$ for $i>0$.
\item $\Hom_A(\Hom_A(M\otimes_BX, A), A)\cong M\otimes_BX$.

Since $$\Hom_A(M\otimes_BX, A)\cong \Hom_B(X, \Hom_A(M, A)).$$
We have seen that $$\Ext^i_B(X, \Hom_A(M, A))=0$$ for $i>0$,
then from Lemma \ref{lemma-defect1} and $X\in \mathcal{G}_B$,  it follows that $$\Hom_B(X, \Hom_A(M, A))\cong \Hom_B(X, B)\otimes_B\Hom_A(M, A)$$
Hence
\begin{eqnarray*}
\Hom_A(\Hom_A(M\otimes_BX, A), A)&\cong&\Hom_A(\Hom_B(X, B)\otimes_B\Hom_A(M, A), A)\\
&\cong& \Hom_B(\Hom_B(X, B), \Hom_A(\Hom_A(M, A), A))\\
&\cong& \Hom_B(\Hom_B(X, B), M)\\
&\cong& M\otimes_B\Hom_B(\Hom_B(X, B), B)\\
&\cong& M\otimes_BX
\end{eqnarray*}

\item $\Ext_A^i(\Hom_A(M\otimes_BX, A), A)=0$ for $i>0.$
We have the following isomorphisms for $i>0$
\begin{eqnarray*}
\Ext_A^i(\Hom_A(M\otimes_BX, A), A) &\cong& \Ext_A^i(\Hom_B(X, \Hom_A(M, A)), A)\\
&\cong& \Ext_A^i(\Hom_B(X, B)\otimes_B\Hom_A(M, A),  A)\\
&\cong& H^i(R\Hom_A(\Hom_B(X, B)\otimes_B^{\LL}\Hom_A(M, A), A)\\
&\cong&H^i(R\Hom_B(\Hom_B(X, B), R\Hom_A(\Hom_A(M, A), A))\\
&\cong& H^i(R\Hom_B(\Hom_B(X, B), M)\\
&\cong& \Ext_B^i(\Hom_B(X, B), M)\\
&\cong& 0
\end{eqnarray*}
where the third isomorphism is the definition of $\Ext^i_A(-, -)$, the forth isomorphism comes from the fact that $$(-\otimes^{\LL}_B\Hom_A(M, A), R\Hom_A(\Hom_A(M, A), -))$$ is an adjoint pair in $\DD^b(B\otimes A^{\op})$ and the seventh isomorphism holds since $X$ is of Gorenstein dimension zero.
\end{enumerate}
Therefore, $M\otimes_BX\in\mathcal{G}_A.$ Similarly, we have $N\otimes_AY\in \mathcal{G}_B$ for any $Y\in \mathcal{G}_A$.

On the other hand, we have known that $\underline{\mathcal{G}_A}$ is a triangulated category with the suspension functor $\Phi_A:=[1]$, which is a quasi-inverse of the syzygy functor $\Omega_A$ (cf. Remark \ref{rem-zero1}).
Therefore we have 
$$M\otimes_B(N\otimes_A-)\cong \Omega_{A}^n\cong [-n]$$ from $\underline{\mathcal{G}_A}$ to $\underline{\mathcal{G}_A}$.
Similarly, we have $$N\otimes_A(M\otimes_B-)\cong \Omega_B^n\cong [-n]$$ from $\underline{\mathcal{G}_B}$ to $\underline{\mathcal{G}_B}.$
So $$M\otimes_B-: \underline{\mathcal{G}_B}\rightarrow \underline{\mathcal{G}_A}$$ is an equivalence of triangulated categories
with an quasi-inverse
 $$N\otimes_A-:\underline{\mathcal{G}_A}\rightarrow \underline{\mathcal{G}_B}.$$
\epf

\begin{prop}\label{prop-derived-sin}
Let $A$ and $B$ be two finite dimensional $k$-algebras such that $\DD^b(A)\cong \DD^b(B)$. Then there exists a pair of bimodules $(_AM_B, _BN_A)$ inducing a singular equivalence of Morita type with level such that they induce an equivalence between $\underline{\mathcal{G}_A}$ and $\underline{\mathcal{G}_B}$,
$$F_M:\underline{\mathcal{G}_B}\rightarrow \underline{\mathcal{G}_A}$$  and $$F_N:\underline{\mathcal{G}_A}\rightarrow \underline{\mathcal{G}_B}.$$
\end{prop}
\pf  Since $\DD^b(A)\cong\DD^b(B)$, we know that there exist two 2-sided tilting complexes:
$$X_{\bullet}: 0\leftarrow X_p\leftarrow X_{p+1}\leftarrow \cdots \leftarrow X_{m-1}\leftarrow M\leftarrow 0$$
and $$Y_{\bullet}:0\leftarrow Y_q\leftarrow Y_{q+1}\leftarrow \cdots\leftarrow Y_{n-1}\leftarrow N\leftarrow 0 $$
such that $X_i$ and $Y_j$ are projective as bimodules for $p\leqslant i\leqslant m-1$ and $q\leqslant j\leqslant n-1$ and $(M, N)$ defines a singular equivalence of Morita type
with level $m+n+1$ (cf. Theorem \ref{thm1} and its proof).
Next we claim that for any $X\in \mathcal{G}_B$, there exists $p\in \calN$ (depending on $X$) such that $$\Omega_A^p(M\otimes_BX)\in \mathcal{G}_A.$$
By assumption that $X\in \mathcal{G}_B$, from Theorem \ref{thm-Ka1} we obtain that 
$$X_p\otimes_BX\leftarrow X_{p+1}\otimes_BX\leftarrow \cdots\leftarrow X_{m-1}\otimes_BX\leftarrow M\otimes_BX\leftarrow 0\in \DD^b(A)_{\fGd}.$$ Since $X_i\otimes_BX$ is projective as an $A$-module for all $i=p, \cdots, m-1$, there exists $p\in \calN$ such that (cf. Proposition 2.10, \cite{Ka}) $$\Omega_A^p(M\otimes_BX)\in \mathcal{G}_A.$$
Therefore we define the functor $F_M$ as follows:
$$F_M(X)\cong \Phi_A^p\Omega_A^p(M\otimes_BX),$$
where $\Phi_A$ is the suspension functor in $\underline{\calG_A}$ (cf. Remark \ref{rem-zero1}).
Note that $F_M$ is a well-defined triangle functor. Similarly, we can define $F_N$.
Moreover, we have
\begin{eqnarray*}
F_NF_M(X)&\cong& \Phi_B^q\Omega_B^q(N\otimes_A\Phi_A^p\Omega_A^p(M\otimes_BX))\\
&\cong&\Phi_B^q(N\otimes_A \Omega_A^q\Phi_A^p\Omega_A^p(M\otimes_BX))\\
&\cong& \Phi_B^q\Phi_B^p\Omega_B^{p+q}(N\otimes_AM\otimes_BX))\\
&\cong& \Phi_B^{p+q}\Omega_B^{p+q}(\Omega_B^n(X))\\
&\cong& \Omega_B^n(X)\\
&\cong& X[-n]
\end{eqnarray*}
Similarly, we have $F_MF_N\cong [-n]$.
Therefore $F_M$ and $F_N$ are quasi-inverse to each other.
\epf

Next we study when a singular equivalence of Morita type with level induces a stable equivalence of Morita type.
We have the following proposition.
\begin{prop}\label{prop-zero3}
Let $(M, N)$ define a singular equivalence of Morita type with level $n$ such that $\Hom_A(M, A)$ and $\Hom_B(N, B)$ are of finite projective dimension as a left $B$-module and as  a left $A$-module respectively.
Suppose that $A\in \calG_{A^e}$ and $B\in \calG_{B^e}$. Then $A$ and $B$ are stably equivalent of Morita type.
\end{prop}

For the proof of Proposition \ref{prop-zero3} above, we need the following two lemmas.
\begin{lemma}\label{lemma-zero3}
Let $X$ be an $A$-$A$-bimodule, $Y$ be a $A$-$B$-bimodule. Then we have the following natural isomorphism
$$\Hom_{A\otimes B^{\op}}(X\otimes_AY, A\otimes B^{\op})\cong \Hom_{A\otimes A^{\op}}(X, \Hom_B(Y, A\otimes B^{\op}).$$
\end{lemma}
\pf First we have
$$\Hom_{A\otimes B^{\op}}(X\otimes_AY, A\otimes B^{\op})=\{f\in \Hom_{B^{\op}}(X\otimes_AY, A\otimes B^{\op}) |  \ f(a(x\otimes y)=af(x\otimes y)
\}$$
From the ``isomorphisme cher \`a Henri Cartan", we have
$$\phi: \Hom_{B^{\op}}(X\otimes_AY, A\otimes B^{\op})\cong \Hom_{A^{\op}}(X, \Hom_{B^{\op}}(Y, A\otimes B^{\op})).$$
We observe that under the isomorphism $\phi$, $f\in \Hom_{B^{\op}}(X\otimes_AY, A\otimes B^{\op})$ such that
$f(a(x\otimes y))=af(x\otimes y)$ if and only if $\phi(f)\in \Hom_{A^{\op}}(X, \Hom_{B^{\op}}(Y, A\otimes B^{\op}))$ such that
$\phi(f)(ax)=a\phi(f)(x)$ for any $a\in A, x\in X$ and $y\in Y$, that is to say,  the isomorphism $\phi$ induces an isomorphism between $$\Hom_{A\otimes B^{\op}}(X\otimes_AY, A\otimes B^{\op})$$ and $$\Hom_{A\otimes A^{\op}}(X, \Hom_B(Y, A\otimes B^{\op}).$$
\epf

\begin{lemma}\label{lemma-zero4}
Let $A$ be a $k$-algebra such that $A\in\calG_{A^e}$, $B$ be a $k$-algebra and $N$ be a $B$-$A$-bimodule such that $\Hom_B(N, B)$ has finite projective dimension
as an $A$-module. Then $N\in \calG_{B\otimes A^{\op}}$.
\end{lemma}
\pf Since $A\in \calG_{A^e}$, there  exists an acyclic complex of projective $A^e$-modules
\begin{eqnarray*}
  \xymatrix{P_{\bullet}: \cdots& P_{-1} \ar[l]_-{d_{-1}}& P_0\ar[l]_-{d_{0}} &P_1\ar[l]_-{d_1} &  P_2\ar[l]_-{d_2} & \cdots\ar[l]}
\end{eqnarray*}
such that
$\ker(d_0)\cong A$ and $\Hom_{A^e}(P_{\bullet}, A^e)$ is acyclic (cf. Remark \ref{rem-zero1}).
Since $N$ is projective as an $A$-module, $N\otimes_A P_{\bullet}$ is a complex of projective $B$-$A$-bimodules and acyclic. Moreover $\ker(N\otimes_A d_0)\cong N.$
Now it is sufficient to show that $\Hom_{B\otimes A^{\op}}(N\otimes_AP_{\bullet}, B\otimes A^{\op})$ is acyclic in order to show that $N\in\calG_{B\otimes A^{\op}}$.
From Lemma \ref{lemma-zero3}, we have the following isomorphisms
\begin{eqnarray*}
\Hom_{B\otimes A^{\op}}(N\otimes_A P_i, B\otimes A^{\op}) &\cong& \Hom_{A^e}(P_i, \Hom_B(N, B\otimes A^{\op}))\\
&\cong&  \Hom_{A^e}(P_i, \Hom_B(N, B)\otimes A^{\op})\\
&\cong& \Hom_{A^e}(P_i, A^e)\otimes_A\Hom_B(N, B).
\end{eqnarray*}
 Next we claim that $\Hom_{A^e}(P_{\bullet}, A^e)\otimes_A\Hom_B(N, B)$ is acyclic. Let us compute the homology group of $\Hom_{A^e}(P_{\bullet}, A^e)\otimes_A\Hom_B(N, B),$ we have the following isomorphism for any  $m\in\calN$ since $\Hom_{A^e}(P_{\bullet}, A^e)$ is an acyclic complex of projective $A^e$-modules,
 \begin{eqnarray*}
   H_i(\Hom_{A^e}(P_{\bullet}, A^e)\otimes_A\Hom_B(N, B))\cong \Tor_{A}^{i+m}(\ker(\Hom_{A^e}(d_{i+m+1}, A^e)), \Hom_B(N, B)).
 \end{eqnarray*}
 When $i+m>\proj.\dim_A (\Hom_B(N, B))$, we have $$\Tor_{A^e}^{i+m}(\ker(\Hom_{A^e}(d_{i+m+1}, A^e)), \Hom_B(N, B))=0,$$ hence $$H_i(\Hom_{A^e}(P_{\bullet}, A^e)\otimes_A\Hom_B(N, B))=0$$ for any $i\in \Z$.
Therefore $\Hom_{B\otimes A^{\op}}(N\otimes_AP_{\bullet}, B\otimes A^{\op})$ is acyclic. So we have shown that $N\in\calG_{B\otimes A^{\op}}.$
\epf

Now we can prove Proposition \ref{prop-zero3}. From the isomorphism $M\otimes_BN\cong \Omega_{A^e}^n(A)$ in $A^e$-$\underline{\modu}$ it follows that
\begin{eqnarray}\label{eqn-wang}
A\cong \Phi_{A^e}^{n}\Omega_{A^e}^n(A)\cong \Phi_{A^e}^{n}(M\otimes_BN)
\end{eqnarray}
in $\underline{\calG_{A^e}}$, where the first isomorphism holds since $A\in \calG_{A^e}$.
From Lemma \ref{lemma-zero4}, we know that $N\in \calG_{B\otimes A^{\op}},$ hence
we have an acyclic resolution of $N$ 
\begin{eqnarray*}
  \xymatrix{Q_{\bullet}: \cdots& Q_{-1} \ar[l]_-{d_{-1}}& Q_0\ar[l]_-{d_{0}} &Q_1\ar[l]_-{d_1} &  Q_2\ar[l]_-{d_2} & \cdots\ar[l]}
\end{eqnarray*}
such that $\ker(d_0)\cong N$ and  $\Hom_{B\otimes A^{\op}}(Q_{\bullet}, B\otimes A^{\op})$ is also acyclic.
Note that $$\ker(d_{-n})\cong \Phi^n_{B\otimes A^{\op}}(N).$$ 
Now we claim that $$M\otimes_B \Phi^n_{B\otimes A^{\op}}(N)\in \calG_{A\otimes A^{\op}}.$$
Since $M$ is projective as a left $A$-module and as a right $B$-module, $M\otimes_B Q_{\bullet}$ is an acyclic 
resolution of $M\otimes_B \Phi^n_{B\otimes A^{\op}}(N)$. It is sufficient to show that 
$\Hom_{A^e}(M\otimes_B Q_{\bullet}, A^e)$ is acyclic. Since we have the following 
isomorphism
\begin{eqnarray*}
\Hom_{A^e}(M\otimes_B Q_i, A^e)&\cong& \Hom_{B\otimes A^{\op}}(Q_i, \Hom_A(M, A^e))\\
&\cong& \Hom_{B\otimes A^{\op}}(Q_i, \Hom_A(M, A)\otimes A^{\op})
\end{eqnarray*}
where the first isomorphism comes from Lemma \ref{prop-zero3}.
Hence $\Hom_{A^e}(M\otimes_B Q_{\bullet}, A^e)$ is acyclic 
since $\Hom_A(M, A)$ is of finite projective dimension as a $B$-module, hence
$$M\otimes_B \Phi^n_{B\otimes A^{\op}}(N)\in \calG_{A^e}.$$
So we have the following isomorphism in $\underline{\calG_{A^e}}$,
 $$\Phi_{A^e}^n(M\otimes_BN)\cong M\otimes_B \Phi^n_{B\otimes A^{\op}}(N),$$
 hence from the isomorphism ($\ref{eqn-wang}$) it follows that we have the following isomorphism in $\underline{\calG_{A^e}}$
 (also in $A^e$-$\underline{\modu}$).
 $$A\cong M\otimes_B \Phi^n_{B\otimes A^{\op}}(N).$$
 
Similarly, we also have
$$B\cong \Phi_{B\otimes A^{\op}}^{n}(N)\otimes_A M$$ in $\underline{\calG_{B^e}}$.
So $(M, \Phi^{n}_{B\otimes A^{\op}}(N))$ defines a stable equivalence of Morita type between $A$ and $B$.
\epf
\subsection{Gorenstein algebras}
In this subsection, we will study singular equivalences of Morita type with
level between Gorenstein algebras. First let us recall some notions on Gorenstein
algebras and finitistic dimensions.
\begin{defn}
  Let $A$ be a Noetherian $k$-algbra. We say that $A$ is Gorenstein if
  the injective dimension of $A$ as a left $A$-module and as a right $A$-module
  is, respectively, finite.
\end{defn}
\begin{defn}
  Let $A$ be a Noetherian $k$-algebra. We define the left (right) finitistic dimension as
  $\sup\{{\proj. \dim P \ |\ \mbox{$P$ is a left (right) $A$-module  of finite projective dimension}}\}.$
\end{defn}
\begin{lemma}[\cite{AuRe}]\label{lemma-fini-AuRe}
  Let $A$ be a $k$-algebras such that $\inj.\dim(_AA)<+\infty$, then $$\inj.\dim(A_A)<+\infty$$ if and only if the right finitistic dimension of $A$ is finite. As the result,  a Gorenstein algebra $A$ has finite left (and right) finitistic dimension.
\end{lemma}

\begin{lemma}\label{prop-Gorenstein-finitistic}
  Let $(M, N)$ define a singular equivalence of Morita type of level $n$ between two $k$-algebras $A$ and $B$. If $A$ has finite left (or right) finitistic dimension, then so does $B$.
\end{lemma}
\pf Denote the left (or right) finitistic dimension of $A$ by $m$. We will prove that the left (or right) finitistic dimension of $B$ is smaller than $m+n$.
Let $X$ be any $B$-module with finite projective dimension, then $M\otimes_BX$ is of finite projective dimension as an $A$-module, hence $\proj.\dim_A(M\otimes_BX)\leqslant m$. Then $$\proj.\dim_B(N\otimes_AM\otimes_BX)\leqslant m.$$
Thus $$\proj.\dim_B(\Omega_B^n(X))\leqslant m,$$ so $$\proj.\dim_B(X)\leqslant m+n.$$
Therefore the finitistic dimension of $B$ is smaller that $m+n$, so it is finite.
\epf

\begin{lemma}\label{lemma-syzygy12}
  Let $A$ be a Noetherian $k$-algebra, let $X$ be an $A$-$A$-bimodule such that $_AX$ is projective as a left $A$-module and $Y$ be a left $A$-module. Then in $A$-$\overline{\modu}$, the stable module category of $A$-$\modu$ modulo injective objects, we have the following isomorphism for any $n\in\calN$.
  $$\Hom_A(\Omega^n_{A^e}(X), Y)\cong \Hom_A(X, \Omega_{A}^{-n}(Y)).$$
\end{lemma}
\pf Take a projective resolution of $X$ in $A^e$-$\modu$,
$$0\leftarrow P_0\leftarrow P_1\leftarrow P_2\leftarrow \cdots,$$
then we have the following exact sequence
$$0\leftarrow X\leftarrow P_0\leftarrow \cdots \leftarrow P_{n-1}\leftarrow \Omega_{A^e}^n(X)\leftarrow 0.$$
Since $X$ is projective as a left $A$-module, we have the following exact sequence in $A$-$\modu$
$$0\rightarrow \Hom_A(X, Y)\rightarrow \Hom_A(P_0, Y)\rightarrow \cdots \rightarrow \Hom_A(\Omega_{A^e}^n(X), Y)\rightarrow 0.$$
Next we claim that $\Hom_A(P_i, Y)$ is an injective $A$-module for all $i=0, \cdots n-1$.
Since $P_i$ is projective as an $A$-$A$-bimodule, there exist a projective $A$-$A$-bimodule $P_i'$ and $m\in\calN$ such that $$P_i\oplus P_i'\cong (A\otimes A^{\op})^m.$$
Hence $$\Hom_A(P_i, Y)\oplus \Hom_A(P_i', Y)\cong \Hom_A(A\otimes A^{\op}, Y)^{\oplus m}.$$
So in order to prove that claim, it is sufficient to show that $\Hom_A(A\otimes A^{\op}, Y)$ is injective as a left $A$-module.
Indeed, we have the following isomorphisms in $A$-$\modu$,
\begin{eqnarray*}
  \Hom_A(A\otimes A^{\op}, Y)&\cong& \Hom_k(A^{\op}, \Hom_A(A, Y))\\
  &\cong& \Hom_k(A^{\op}, k)\otimes \Hom_A(A, Y)
\end{eqnarray*}
where the first isomorphism is the ``isomorphisme cher \`a Henri Cartan", hence we know
that $\Hom_A(A\otimes A^{\op}, Y)$ is an injective $A$-module, so we have shown the claim.
 We have the following isomorphism in $A$-$\overline{\modu}$ since the co-syzygies are independent of the injective resolution of $\Hom_A(X, Y)$ in $A$-$\overline{\modu}$.
\begin{eqnarray}\label{eqn12}
  \Omega_A^{-n}(\Hom_A(X, Y))\cong \Hom_A(\Omega_{A^e}^n(X), Y).
\end{eqnarray}

On the other hand, take an injective resolution of $Y$ in $A$-$\modu$,
$$0\rightarrow I_0\rightarrow I_1\rightarrow \cdots, $$
then we have the following exact sequence
$$0\rightarrow Y\rightarrow I_0\rightarrow \cdots \rightarrow I_{n-1}\rightarrow \Omega_A^{-n}(Y)\rightarrow 0.$$
Since $X$ is projective as a left $A$-module, we have the following exact sequence
$$0\rightarrow \Hom_A(X, Y)\rightarrow \cdots\rightarrow \Hom_A(X, I_{n-1})\rightarrow \Hom_A(X, \Omega_A^{-n}(Y))\rightarrow 0.$$
Since $I_i$ is injective, $\Hom_A(X, I_i)$ is injective for $i=0, \cdots n-1$, hence in $A$-$\overline{\modu}$ we have the following isomorphism \begin{eqnarray}\label{eqn13}
  \Omega_A^{-n}(\Hom_A(X, Y))\cong \Hom_A(X, \Omega_A^{-n}(Y)).
\end{eqnarray}
Hence from the isomorphisms (\ref{eqn12}) and (\ref{eqn13}), we have the isomorphism in $A$-$\underline{\modu}$,
  $$\Hom_A(X, \Omega_A^{-n}(Y))\cong\Hom_A(\Omega_{A^e}^n(X), Y).$$
  \epf

\begin{prop}\label{prop-Gorenstein-duality}
 Let $A$ and $B$ be two Noetherian $k$-algebras, let $(M, N)$ define a singular equivalence of Morita type of level $n$ between $A$ and $B$. Suppose that $\Hom_B(N, B)$ has finite projective dimension as $A$-module. If $A$ is Gorenstein, so is $B$.
\end{prop}
\pf  Let $X$ be a projective $B$-module, then $\Hom_B(N, X)$ has finite projective dimension as an $A$-module since $\Hom_B(N, B)$ has finite projective dimension as an $A$-module. Since $A$ is a Gorenstein algebra, $\Hom_B(N, X)$ has finite injective dimension. Since $M$ is projective as left $A$-module and as a right $B$-module, $\Hom_A(M, -)$ is exact and sends injective modules to injective modules. It follows that $$\Hom_A(M, \Hom_B(N, X))$$ has finite injective dimension as $B$-module.
On the other hand, we have the following isomorphisms in $B$-$\overline{\modu}$, the stable module category of $B$-$\modu$ modulo injective objects.
\begin{eqnarray*}
  \Hom_A(M, \Hom_B(N, X))&\cong& \Hom_B(N\otimes_AM, X)\\
  &\cong& \Hom_B(\Omega^n_{B^e}(B), X)\\
  &\cong & \Hom_B(B, \Omega_B^{-n}(X))\\
  &\cong& \Omega_B^{-n}(X)
\end{eqnarray*}
where the third isomorphism is because of Lemma \ref{lemma-syzygy12}.
Hence $\Omega_B^{-n}(X)$ has finite injective dimension, thus $X$ has finite injective dimension. So we have proved that each projective $B$-module has finite injective dimension, it follows that $B$ is a left Gorenstein algebra. From Lemma \ref{prop-Gorenstein-finitistic}, we know that the left (or right) finitistic dimension of $B$ is finite since the left (or right) finitistic dimension of $A$ is finite. Thus it follows from Lemma \ref{lemma-fini-AuRe} that $B$ is Gorenstein.
\epf

\begin{cor}
  Let $A$ and $B$ be two Noetherian $k$-algebras, let $(M, N)$ define a singular equivalence of Morita type of level $n$ between $A$ and $B$. Suppose that $\Hom_B(N, B)$ is projective as an $A$-module. If $A$ is Gorenstein with $\inj.\dim(A)=m$, then $B$ is Gorenstein with $\inj.\dim(B)\leqslant m+n$.
\end{cor}
\pf
Since $\Hom_B(N, B)$ is projective as an $A$-module and $A$ is Gorenstein, we have that $\inj. \dim(\Hom_B(N, B))\leqslant m$. From Proposition \ref{prop-Gorenstein-duality},
we have the following isomorphism
\begin{eqnarray*}
  \Hom_A(M, \Hom_B(N, X))
  \cong \Omega_B^{-n}(X)
\end{eqnarray*}
Since $\Hom_A(M, -)$ is exact and sends injective $A$-modules to injective $B$-modules, we have
$$\inj.\dim( \Hom_A(M, \Hom_B(N, X)))\leqslant m,$$
hence $$\inj.\dim(\Omega^{-n}(B))\leqslant m,$$ so $\inj.\dim(B)\leqslant m+n$.
\epf

\appendix
\section{Singular equivalence of Morita type in the Sense of Chen and Sun}

%
%
%
%
%

Recall the definition of singular equivalence of Morita type (cf. Definition \ref{def-chen}) defined by Chen and Sun. It is a strong version of singular equivalence of Morita type with level, which means, a singular equivalence of Morita type in the sense of Chen and Sun induces a singular equivalence of Morita type with level.
In general, the self-injectivity of algebras may not be preserved under singular equivalence of Morita type with level (cf. Remark \ref{rem-expl-self}). Next we will prove that it is preserved under singular equivalence of Morita type.
\begin{prop}\label{prop-singular-Sun}
  Let $A$ and $B$ be two finite dimensional algebras over a field $k$, $(M, N)$ defines a singular equivalence of Morita type (in the sense of Chen and Sun) between $A$ and $B$. If $A$ is self-injective, then so is $B$.
\end{prop}
Before the proof of Proposition \ref{prop-singular-Sun},
we recall some notions about singularly stable category (cf. \cite{ZhZi}).
\begin{defn}[Definition 1.2, \cite{ZhZi}]
  Let $A$ be a finite dimensional $k$-algebra, We denote by $\mathcal{P}^{<\infty}(A)$ the full
  subcategory of $A$-$\modu$ consisting of $A$-modules with finite projective dimension. Then we
  define the singularly stable category as the quotient category of $A$-$\modu$ by $\calP^{<\infty}(A)$, denoted by $A$-$\underline{\modu}_{\calP^{<\infty}}$.
\end{defn}
\begin{rem}
There is a natural functor $$L_A: \mbox{$A$-$\underline{\modu}_{\calP^{<\infty}}$}\rightarrow \DD_{\sg}(A).$$ Moreover, we have $$[-1]\circ L_A\cong L_A\circ \Omega_A.$$
  We also remark that a singular equivalence of Morita type induces an equivalence between singularly
  stable categories. If $A$ is self-injective, then $\mbox{$A$-$\underline{\modu}_{\calP^{<\infty}}$}\cong \mbox{$A$-$\underline{\modu}$}$.
\end{rem}
\begin{lemma}\label{lemma-chen1}
  Let $(M, N)$ define a singular equivalence of Morita type (in the sense of Chen and Sun)
  between two finite dimensional $k$-algebras $A$ and $B$. If $A$ is self-injective, then the natural
  functor  $$L_B: \mbox{$B$-$\underline{\modu}_{\calP^{<\infty}}$}\rightarrow \DD_{\sg}(B)$$
  is an equivalence. Moreover, $$\Omega_B: \mbox{$B$-$\underline{\modu}_{\calP^{<\infty}}$}\rightarrow \mbox{$B$-$\underline{\modu}_{\calP^{<\infty}}$}$$ is an equivalence.
\end{lemma}
\pf Since $(M, N)$ define a singular equivalence of Morita type between $A$ and $B$, we have the following commutative diagram:
\begin{eqnarray*}
  \xymatrix{
  \DD_{\sg}(B)\ar[r]_{\cong}^{M\otimes_B-} & \DD_{\sg}(A)\\
  \mbox{$B$-$\underline{\modu}_{\calP^{<\infty}}$}\ar[u]^{L_B}\ar[r]_{\cong}^{M\otimes_B-} &
  \mbox{$A$-$\underline{\modu}_{\calP^{<\infty}}$}\ar[u]_{L_A}^{\cong}
  }
\end{eqnarray*}
Since $A$ is self-injective, $L_A$ is an equivalence, hence $L_B$ is also an equivalence.
From $[-1]\circ L_B\cong L_B\circ \Omega_B,$ it follows that
$$\Omega_B: \mbox{$B$-$\underline{\modu}_{\calP^{<\infty}}$}\rightarrow \mbox{$B$-$\underline{\modu}_{\calP^{<\infty}}$}$$ is an equivalence.
\epf
\begin{lemma}\label{lemma-chen2}
  Let $(M, N)$ define a singular equivalence of Morita type between $A$ and $B$ and $A$
  be self-injective, then every injective $B$-module is of finite projective dimension,
  moreover, $M\otimes_B-$ sends injective $B$-modules to injective $A$-modules.
\end{lemma}
\pf Suppose that $I$ is an indecomposable injective $B$-module with infinite projective dimension, then
$I\neq 0$ in $B$-$\underline{\modu}_{\calP^{<\infty}}$. From Lemma \ref{lemma-chen1}, we know
that $$\Omega_B: \mbox{$B$-$\underline{\modu}_{\calP^{<\infty}}$}\rightarrow \mbox{$B$-$\underline{\modu}_{\calP^{<\infty}}$}$$ is an equivalence. Hence there exists a
$B$-module $X$ with infinite projective dimension such that $$\Omega_B(X)\cong I$$ in
$B$-$\underline{\modu}_{\calP^{<\infty}}$, that is, there exist a $B$-module of finite
projective dimension $P_1$ such that $$\Omega_B(X)\cong I\oplus P_1$$ in $B$-$\underline{\modu}$ since $I$ is indecomposable.
Hence we have an exact sequence
\begin{eqnarray*}
  0\rightarrow I\oplus P_1\rightarrow P_0\rightarrow X\rightarrow 0
\end{eqnarray*}
where $P_0$ is projective. So it follows that $I$ is a submodule of $P_0$. Since $I$ is injective, $I$ is a direct summand of $P_0$, hence $I$ is projective, which is a contradiction
with that $I$ is of infinite projective dimension.
Therefore each injective $B$-module has finite projective dimension. Suppose that $M\otimes_BI$ is not an injective $A$-module, that is $M\otimes_BI\neq 0$ in $A$-$\underline{\modu}$. Since
$$N\otimes_A-: \mbox{$A$-$\underline{\modu}$}\rightarrow \mbox{$B$-$\underline{\modu}_{\calP^{<\infty}}$}$$ is an equivalence,
$N\otimes_AM\otimes_BI\neq 0$ in $B$-$\underline{\modu}_{\calP^{<\infty}}$, since $N\otimes_AM\otimes_B I\cong I\oplus P\otimes_BI$,
$$I\oplus P\otimes_B I\neq 0$$ in $B$-$\underline{\modu}_{\calP^{<\infty}}.$

 Since $P$ is
projective as a left $B$-module and as a right $B$-module and $I$ is of finite projective dimension,  we have $\proj.\dim P\otimes_BI <\proj.\dim I<\infty$. Hence $$I\oplus P\otimes_B I=0$$ in $B$-$\underline{\modu}_{\calP^{<\infty}}.$  Contradiction! Therefore $M\otimes_BI$ is an injective $A$-module for any injective $B$-module $I$.
\epf

{\it Proof of Proposition \ref{prop-singular-Sun}.\hspace{2ex}}\ \
Let $I$ be an injective $B$-module, from Lemma \ref{lemma-chen2}, we know that $M\otimes_BI$ is projective (injective),  thus $$N\otimes_AM\otimes_BI\cong I\oplus Q\otimes_BI$$ are projective, hence $I$
is projective. So we have proved that every injective $B$-module is projective, therefore $B$ is self-injective.
\epf

\bibliographystyle{plain}

\end{document}